\documentclass[12pt]{amsart}
\usepackage{amscd,amsmath,amssymb,amsfonts}
\usepackage[all]{xy}
\theoremstyle{plain}
\newtheorem{thm}{Theorem}
\newtheorem{lem}[thm]{Lemma}

\newtheorem{prop}[thm]{Proposition}

\theoremstyle{definition}
\newtheorem{defn}[thm]{Definition}

\newtheorem{rmk}[thm]{Remark}
\newtheorem{rmks}[thm]{Remarks}

\newtheorem{nota}[thm]{Notation}

\numberwithin{thm}{section}
\numberwithin{equation}{section}

\newcommand{\ga}[2]{\begin{gather}\label{#1}#2 \end{gather}}

\newcommand{\surj}{\twoheadrightarrow}
\newcommand{\inj}{\hookrightarrow}

\newcommand{\Hom}{{\rm Hom}}
\newcommand{\im}{{\rm im}}

\newcommand{\sC}{{\mathcal C}}
\newcommand{\sD}{{\mathcal D}}
\newcommand{\sE}{{\mathcal E}}
\newcommand{\sF}{{\mathcal F}}
\newcommand{\sG}{{\mathcal G}}

\newcommand{\sN}{{\mathcal N}}
\newcommand{\sO}{{\mathcal O}}

\newcommand{\sQ}{{\mathcal Q}}

\newcommand{\sS}{{\mathcal S}}
\newcommand{\sT}{{\mathcal T}}

\renewcommand{\P}{{\mathbb P}}

\newcommand{\et}{{\acute{e}t}}
\newcommand{\Rep}{{\rm Rep\hspace{0.1ex}}}

\begin{document}

\title{On Nori's fundamental group scheme}
\author{H\'el\`ene Esnault}
\address{
Universit\"at Duisburg-Essen, Mathematik, 45117 Essen, Germany}
\email{esnault@uni-due.de}
\author{Ph\`ung  H\^o Hai}
\address{
Universit\"at Duisburg-Essen, Mathematik, 45117 Essen, Germany
and Institute of Mathematics, Hanoi, Vietnam}
\email{hai.phung@uni-duisburg-essen.de}
\author{Xiaotao Sun}
\address{Chinese Academy of Mathematics and Systems Science, Beijing, P. R. of China}
\email{xsun@math.ac.cn}
\address{}
\date{June 30, 2006}
\thanks{Partially supported by  the DFG Leibniz Preis,  the DFG Heisenberg program and the grant NFSC 10025103}
\begin{abstract}
The aim of this note is to give two structure theorems on Nori's fundamental group scheme of a proper connected variety defined over a perfect field
and endowed with a rational point.
\end{abstract}
\maketitle
\begin{quote}

\end{quote}

\section{Introduction}
For a proper connected reduced scheme $X$ defined over a perfect field $k$ endowed with a rational point $x\in X(k)$,
Nori defined in \cite{N1} and \cite{N2} a fundamental group scheme $\pi^N(X,x)$ over $k$.
It is Tannaka dual to the $k$-linear abelian rigid tensor category $\sC^N(X)$ of {\it Nori finite} bundles,
that is bundles which are trivializable over a principal bundle $\pi: Y\to X$ under a finite group scheme.
The rational point $x$ endows $\sC^N(X)$ with a fiber functor $V\mapsto V|_x$ with values in the category of
finite dimensional vector spaces over $k$. This makes $\sC^N(X)$  a Tannaka category, thus by Tannaka duality (\cite[Theorem~2.11]{DeMil}),
the fiber functor establishes an equivalence between $\sC^N(X)$ and the representation category ${\rm Rep}(\pi^N(X,x))$ of an affine group scheme
$\pi^N(X,x)$ (see Section 2 for an account of Nori's construction), that is a prosystem of affine algebraic $k$-group schemes, which turn out to be finite group schemes. 
The purpose of this note is to study the structure of this Tannaka group scheme.

To this aim, we define two full tensor subcategories $\sC^\et(X)$ and $\sC^F(X)$.
The objects of the first one are {\it \'etale finite} bundles, that is bundles for which the finite group scheme is \'etale,
while the objects of the second one are $F$-{\it finite} bundles, that is bundles for which the group scheme is local.
As sub-Tannaka categories they are
representation categories of Tannaka group schemes $\pi^\et(X,x)$ and $\pi^F(X,x)$.  Our first main theorem (see Theorem \ref{thm4.1})
asserts that the natural homorphism of
$k$-group schemes
\ga{1.1}{\pi^N(X,x)\to \pi^\et(X,x)\times \pi^F(X,x)}
is faithfully flat, so in particular surjective.
To have a feeling for the meaning of the statement, it is useful to compare $\pi^\et(X,x)$ with the more familiar fundamental group $\pi_1(X, \bar{x})$
defined by Grothendieck in
(\cite[Expos\'e~5]{Groth}), where $\bar{x}$ is a geometric point above $x$. Grothendieck's fundamental group is
an abstract group, which is a pro-system of finite abstract groups. One has
\ga{1.2}{\pi^\et(X, x)(\bar{k})\cong \pi_1(X\times_k \bar{k},\bar{x})}
(see Remarks \ref{rmks3.2} for a detailed discussion), thus the  \'etale piece of Nori's group scheme takes into
account only the geometric fundamental group and ignores somehow  arithmetics. On the other hand, $\pi^F(X,x)$ reflects
the purely inseparable covers of $X$. That $k$ is perfect guarantees that inseparable covers come only from geometry,
and not from the ground field.

However \eqref{1.1} is not injective, and we give in Remark \ref{rmk4.3} an example  based on Raynaud's work \cite{R} on coverings of
curves producing a new ordinary part in the Jacobian.

The central theorem of our note is the determination by its objects and morphisms
of a
$k$-linear abelian rigid tensor category $\sE$, which is
equivalent to the representation category of ${\rm
Ker}(\eqref{1.1})$ (see Definition \ref{defn5.4} for the
construction and Proposition \ref{prop5.5} and Theorem
\ref{thm5.7} to see that it computes what one wishes). 
This is the most delicate part of the construction. If $S$ is a finite subcategory of $\sC^N(X)$ with \'etale finite Tannaka group scheme $\pi(X, S,x)$, then the total space $X_S$  of the $\pi(X,S,x)$-principal bundle $\pi_S: X_S\to X$  which trivializes all the objects of $S$ has the same property as $X$. It is proper, reduced and connected. However, if $S$ is finite but $\pi(X, S,x)$ is not \'etale, then Nori shows that $X_S$ is still proper connected, but it is not the case that $X_S$ is still reduced. We give a concrete example in Remark \ref{rmk2.4}, 2), which is due to P. Deligne. However, in order to describe $\sE$, we need in some sense an extension of Nori's theory to those non-reduced covers. We define on each such $X_S$ a full subcategory  $\sF(X_S)$ of the category of coherent sheaves, the objects of  which have the property that their push down on $X$ lies in $\sC^N(X)$ (see Definition  \ref{defn2.5}). We show that indeed those coherent sheaves have to be vector bundles (Proposition \ref{prop2.8}), so in a sense, even if the scheme $X_S$ might be bad,  objects which push down to Nori's bundles on $X$
are still good. In particular, $\sC^N(X_S)=\sF(X_S)$ if $\pi(X,S,x)$ is \'etale (Theorem \ref{thm2.10}), so the definition generalizes slightly Nori's theory. 
For given finite subcategories $S$ and $T$ of $\sC^N(X)$, with   $\pi(X,S,x)$ \'etale and $\pi(X,T,x)$ local, one defines in Definition 
\ref{defn5.2} 
 a full subcategory $\sE(X_{S \cup T})\subset \sF(X_{S\cup T})$ on the total space $X_{S\cup T}$ of the principal bundle $\pi_{S\cup T}: X_{S\cup T}\to X$ consisting 
 of those bundles $V$, the push down of which on $X_S$ is $F$-finite. Now the objects of 
$\sE$ are pairs $(X_{S\cup T}, V)$ for $V$ an object in $\sE(X_{S\cup T})$.  
 Morphisms are 
subtle as they do take into account the whole inductive system of
such $T'\subset \sC^F(X)$.

We now describe our method of proof. We proceed in three steps. To see that
$\pi^N(X,x)\to \pi^\et(X,s)$ is surjective is very easy, thus we
consider the kernel $L(X,x)$ and determine its representation
category in Section 3. The computation is based on two results.
The first one of geometric nature asserts that sections of an
$F$-finite bundle can be computed on any principal bundle $X_S\to
X$ with finite \'etale group scheme (see Proposition
\ref{prop3.4}). The second one is the key to the categorial work
and comes from  \cite[Theorem 5.8]{EP}. (For the reader's convenience, we give a short account of the categorial statement in Appendix A).
It gives a criterion for a
$k$-linear abelian rigid category $\sQ$, endowed with a tensor
functor $q: \sT\to \sQ$, to be the quotient category of the full
embedding of $k$-linear  rigid tensor categories
$\sS\xrightarrow{\iota} \sT$, in the sense that if the three
categories are endowed with compatible fiber functors to the
category of $k$-vector spaces, then the Tannaka group schemes
$G(?)$ are inserted in the exact sequence 
\ga{}{1\to
G(\sQ)\xrightarrow{q^*} G(\sT)\xrightarrow{\iota^*} G(\sS)\to
1.\notag} 
The objects of ${\rm Rep}(L(X,x))$ are pairs $(X_S, V)$
where $X_S\to X$ is a principal bundle under an \'etale finite
group scheme and $V$ is a $F$-finite bundle on $X_S$. Morphisms are
defined naturally via Proposition \ref{prop3.4}. The second step
consists in showing surjectivity $L(X,x)\to \pi^F(X,x)$ (see
Theorem \eqref{thm4.1}). To this aim, one needs a strengthening of
Proposition \ref{prop3.4} (see Theorem \ref{thm4.2}) which asserts
that not only sections can be computed on finite \'etale principal
bundles, but also all subbundles of an $F$-finite bundle. Finally
the last step consists in showing that the category $\sE$
constructed in Section 5 is indeed the right quotient category,
for which we use the already mentioned criterion  \cite[Theorem
5.8]{EP}.\\[.2cm]
{\it Acknowledgements:} 
Pierre Deligne sent us his enlightening example which we reproduced in Remark \ref{rmk2.4}, 2).  It allowed us to correct the  main  definition of our category $\sE$ (Section 5)) which was wrongly stated in the first version of this article. We profoundly thank him for  his interest, his encouragement and his help. We also warmly thank  
 Michel Raynaud for answering all our questions on his work on theta characteristics on curves and
on their Jacobians.
The third named author was supported by the Leibniz Preis during the preparation of this work.
\section{ Nori's category}
Let $X$ be a proper reduced scheme defined over a field $k$, which
is {\it connected}
 in the sense that $H^0(X, \sO_X)=k$. We assume  that $k$ has either characteristic 0 or
 that $k$ is perfect of characteristic $p>0$. 
In \cite{N1}, \cite{N2}, Nori defines a category of ``essentially
finite vector bundles" which we recall now.

A vector bundle on $X$ is called \emph{semi-stable} if it is
semi-stable of degree 0 while restricted to each  proper curve in $X$. This is a
full  subcategory of the category $\text{Qcoh}(X)$ of
quasi-coherent sheaves on $X$, is abelian \cite[Lemma 3.6]{N1} and will be
denoted by $\sS(X)$. A vector bundle $V$ on $X$ is called
\emph{finite} if there are polynomials $f\neq g$ whose
coefficients are non-negative integers such that $f(V)$ and $g(V)$
are isomorphic. Nori proves that finite bundles are
semi-stable \cite[Corollary 3.5]{N1} and defines $\sC^N(X)$ to be
the full subcategory of $\sS(X)$ generated by all of the finite
bundles on $X$ in the following sense.
\begin{defn}\label{defn2.1}
The category $\sC^N(X)$ is defined as a full subcategory of the
category  $\sS(X)$. Its objects  are semistable bundles $V$ on
$X$ which have the property that there is a finite collection of
finite-bundles $W_1,\,...,\,W_t$ and semistable bundles
\ga{2.1}{V_1\subset V_2 \subset \bigoplus^t_{i=1}W_i} such that
$V$ is isomorphic to $V_2/V_1$. These are the so called {\it
essentially finite} bundles on $X$.
\end{defn}

If $X$ is a scheme of finite type over $k$, recall that for an affine group scheme $G$ over $k$, $j:P\to X$ is
said to be \emph{a principal} $G$-\emph{bundle on} $X$ 
if 
\begin{itemize}
\item[(1)] $j$ is a faithfully flat affine morphism 
\item[(2)] $\phi:P\times
G\to P$ defines an action of $G$ on $P$ such that
$j\circ\phi=j\circ p_1$ 
\item[ (3)] $(p_1,\phi):P\times G\to P\times_X
P$ is an isomorphism. \end{itemize}
Given such a principal $G$-bundle $P$, one defines a tensor functor \ga{2.2}{\eta_P:\Rep(G)\to
\text{Qcoh}(X),\quad E\mapsto P\times^G E. } 
Conversely, any tensor
functor $\eta:\Rep(G)\to \text{Qcoh}(X)$ yields a principal
$G$-bundle on $X$, which is the spectrum of the $\sO_X$-algebra
$\eta(k[G])$, where $k[G]$ denotes the affine algebra on $G$. See
\cite[Chapter I]{N1} for details.

It is proved that $\sC^N(X)$ is an abelian, rigid $k$-linear
tensor category \cite[Proposition 3.7]{N1}. We sssume that $X$ has a
$k$-rational point $x$, and we fix it throughout the
article. Then one has a fiber functor 
\ga{2.3}{|_x: \sC^N(X)\to
\text{Vect}_k, \ V\mapsto V|_x:=V\otimes_{\sO_X}k} with values in
the category of finite dimensional vector spaces, which defines  a
Tannaka structure on $\sC^N(X)$. We denote by $\pi^N(X,x)$ the
corresponding Tannaka group scheme over $k$. By Tannaka duality
 (\cite[Theorem~2.11]{DeMil}),  one has an equivalence of categories
\ga{2.4}{\sC^N(X)\xrightarrow{|_x \ \ \cong} {\rm
Rep}(\pi^N(X,x)).} We denote by $\eta$ the inverse functor
\ga{2.5}{\eta: {\rm Rep}(\pi^N(X,x))\to  \sC^N(X).} It determines
uniquely a (universal) principal $\pi^N(X,x)$-bundle
$\tilde\pi:\widetilde {X}\to X$ and a $k$-rational point $\tilde
x\in \widetilde X(k)$ over $x\in X(k)$ such that the functor
$\eta$ sends any representation $E$ of $\pi^N(X,x)$ to the
associated vector bundle $\widetilde X\times^{\pi^N(X,x)}E$.

Nori's construction is more precise. Let $S$ be a full tensor
subcategory of $\sC^N(X)$ generated by finitely many objects and
denote its
 Tannaka group by $\pi(X,S,x)$. Then
$\pi(X,S,x)$ is a finite group scheme and
\ga{2.6}{ \pi^N(X,x)=\varprojlim_{S-\text{fin.gen.}} \pi(X,S,x),\quad \pi^N(X,x) \surj \pi(X,S,x).}
Furthermore, the functor $\eta$ applied to the affine algebra $k[\pi(X,S,x)]$
of the group scheme $\pi(X,S,x)$ yields a $\pi(X,S,x)$-principal bundle
\ga{2.7}{\pi_S: X_S\to X, \ {\rm so} \ \widetilde{X}=\varprojlim_{S-\text{fin.gen.}} X_S, \ \tilde{\pi}=\varprojlim_{S-\text{fin.gen.}} \pi_S.}
The scheme $X_S$ is connected in the sense that
\ga{2.8}{H^0(X_S, \sO_{X_S})=k,}
 (\cite[Chapter~II, Proposition~3]{N2}).
Functoriality of $\eta$ with respect to morphisms $Y\to X$, while
applied to $x\to X$, yields that 
\ga{2.9}{x\times_X X_S \to x \
{\rm is \ the \ trivial} \ \pi(X,S,x)-{\rm principal \ bundle}.}
Thus the unit element $x_S$ of this trivial principal bundle,
viewed as a closed point of  $X_S$, is a rational point
\ga{2.10}{x_S\in X_S(k).} Furthermore, $\pi_S$ is universal in the
following sense:
\ga{2.11}{V\in {\rm Obj}(S)\Longleftrightarrow
\pi_S^*(V) \ {\rm trivial}.} 
See also Remark \ref{rmk2.2} below. (We refer to \cite[ Chapters~I,
II]{N2} for the exposition above).

\begin{rmk}\label{rmk2.2} Recall that we assume  $X(k)\neq \emptyset$.
 The category $\sC^N(X)$ can be defined
as a full subcategory of $\sS(X)$ whose objects have the property
that there is a finite group scheme $G$ and a principal $G$-bundle
$\pi: Y\to X$, so that $\pi^*(V)$ is trivializable as an algebraic
bundle. Indeed, if $V$ has this property, then there exists an
injective map $V\inj\pi_*\pi^*V\cong \pi_*\sO_Y^{\oplus d}$ where
$d$ is the rank of $V$. Since $\pi_*\sO_Y^{\oplus d}$ is finite
(\cite[Proposition 3.8]{N1}), $V$ is essentially finite.
Conversely, by Nori's work, any essentially finite bundle
satisfies the property. From now on, we call a bundle with such
property a \emph{ Nori finite} bundle.
\end{rmk}

For a finite bundle $V$, denote by $\langle V\rangle$ the full tensor subcategory generated by $V$.
\begin{defn}\label{defn2.3}
An {\it \'etale finite} bundle is a Nori finite bundle for which $\pi(X,\langle V\rangle,x)$ is \'etale (equivalently is smooth).
If $k$ has characteristic $p>0$, an $F$-{\it finite} bundle is a Nori finite bundle
for which $\pi(X,\langle V\rangle,x)$ is local.
We denote by $\sC^{\acute{e}t}(X)$, resp. $\sC^F(X)$,  the full
tensor subcategory of $\sC^N(X)$ of \'etale, resp. $F$-finite bundles.
\end{defn}
The categories $\sC^\et(X)$ and $\sC^F(X)$ are both tensor subcategories,
thus via the fiber functor at $x$ they yields Tannaka group $\pi^{\acute{e}t}(X,x)$
and $\pi^F(X,x)$, respectively. Furthermore, one has
\ga{2.12}{\sC^{\acute{e}t}(X)\cap \sC^F(X)= \{{\rm trivial \ objects}\}}
where $\{$trivial objects$\}$ means the full subcategory of $\sC^N(X)$
consisting of trivializable bundles, that is isomorphic as bundles to a  direct sum of $\sO_X$.

\begin{rmk} \label{rmk2.4}1)
It is shown in \cite[Section 2]{MS} that a Nori finite bundle $V$  is $F$-finite
if and only if there is a natural number $N>0$ so that $(F_{{\rm abs}}^N)^*(V)$
is trivial, where $F_{{\rm abs}}$ is the absolute Frobenius.\\[.2cm]
2) If 
$S\subset \sC^N(X)$ is a finite subcategory, with $\pi(X, S,x)$ \'etale, then 
$X_S$ is still proper, reduced and connected. However, if $\pi(X, S,x)$ is finite but not \'etale, $X_S$ is still proper and connected \cite[Chapter II, Proposition 3]{N2}, but not necessarily reduced. 
Indeed there are principal bundles $Y\to X$ under a finite local group scheme, so that the total space $Y$ is not reduced.
 We reproduce here a example due to P. Deligne.  Let $k$ be algebraically closed of characteristic $p>0$, and let  $X\subset \P^2$ 
be the union of a smooth conic $X'$ and a tangent line $X''$. Thus $X'\cap X''$ is isomorphic to ${\rm Spec} \ k[\epsilon]/(\epsilon^2)$ as a $k$-scheme. One constructs $\pi: Y\to X$ by gluing the trivial $\mu_p$-torsors $X'\times_k \mu_p$ to $X''\times_k \mu_p$ along a non-constant section of $
{\rm Spec} \ k[\epsilon]/(\epsilon^2)\times_k \mu_p \to {\rm Spec} \ k[\epsilon]/(\epsilon^2)$. For example, one may take the non-constant section ${\rm Spec} \ k[\epsilon]/(\epsilon^2)
\to {\rm Spec} \ k[\epsilon]/(\epsilon^2) \times_k \mu_p$ defined by $
k[\xi, \epsilon]/(\xi^p-1, \epsilon^2)\to k[\epsilon]/(\epsilon^2), \xi \mapsto 
1+\epsilon$. Then $Y$ is projective, non-reduced, and yet fulfills the condition $H^0(Y, \sO_Y)=k$.

\end{rmk}

If $X_S$ is not reduced, there is no good notion of semistable
vector bundles on $X_S$. However, for later use of this article,
we need a category $\sF(X_S)$ similar to $\sC^N(X)$ on the
principal $\pi(X, S,x)$-bundle $\pi_S:X_S\to X$ where $S\subset
\sC^N(X)$ is a finitely generated full tensor subcategory. 
\begin{defn}\label{defn2.5} Let $S\subset \sC^N(X)$ be  a finitely generated full tensor subcategory. Then
$\sF(X_S)\subset \text{Qcoh}(X_S)$ is defined to be the full
subcategory of $\text{Qcoh}(X_S)$, the objects of which are quasi-coherent
sheaves $V$ on $X_S$ such that $(\pi_S)_*V\in\sC^N(X)$.
\end{defn}

It is clear that $\sF(X_S)$ is an abelian category. We will show below that objects of $\sF(X_S)$
come from representations of finite group schemes. In particular,
they are vector bundles on $X_S$, and in fact $\sF(X_S)$ is an abelian rigid $k$-linear tensor category. When $\pi(X,S,x)$ is reduced, we
will show that $\sF(X_S)$ coincides with $\sC^N(X_S)$.

Let $S\subset S'\subset \sC^N(X)$ be finitely generated 
full tensor subcategories. Then one has the following commutative diagram
\ga{2.13}{\xymatrix{ X_{S'}\ar[rr]^{\pi_{S',S}}\ar[dr]_{\pi_{S'}}& &X_S\ar[dl]^{\pi_S}\\
&X }} 
Further one has a faithfully flat homomorphism
\ga{2.14}{G_{S'}:=\pi(X,S',x)\to \pi(X,S,x)=:G_S}

\begin{lem}\label{lem2.6}The morphism $\pi_{S',S}:X_{S'}\to X_S$ is a principal
bundle under a group $G_{S',S}$ which is the kernel of the homomorphism
\eqref{2.14}.
\end{lem}
\begin{proof} Since $X_S$ is the principal bundle induced from
$X_{S'}$ \cite[Proposition 3.11]{N1} via the homomorphism
\eqref{2.14}, it is enough to check that the homomorphism
$G_{S'}\to G_S$ in \eqref{2.14} makes $G_{S'}$ into a principal
$G_{S',S}$-bundle over the scheme $G_S$, which we omit.
\end{proof}

The principal bundle $\pi_{S',S}: X_{S'}\to X_S$ yields a tensor
functor 
\ga{2.15}{\eta_{S',S}:\Rep(G_{S',S})\to
\text{Qcoh}(X_S),\quad \eta_{S',S}(V):=X_{S'}\times^{G_{S',S}} V.}

\begin{lem}\label{lem2.7}The functor $\eta_{S',S}$ in \eqref{2.15} is fully faithful.
Consequently $G_{S',S}$ is isomorphic to the Tannaka group of the
category $\im(\eta_{S',S})$.
\end{lem}
\begin{proof}
It is enough to check that $\eta:=\eta_{S',S}$ is full, i.e., any
morphism $\eta(V)\to \eta(W)$ in $\text{Qcoh}(X_S)$ is induced by
a morphism $V\to W$ in $\text{Rep}(G_{S',S})$.  This is equivalent
to showing $H^0(X_S,\eta(V))\cong V^{G_{S',S}}$ for any $V\in
\text{Rep}(G_{S',S})$. Recall that
$\eta(V):=X_{S'}\times^{G_{S',S}}V$. Thus
\ga{2.16}{H^0(X_S,\eta(V))\cong H^0(X_{S'},\sO_{X_{S'}}\otimes_k
V)^{G_{S',S}}\cong V^{G_{S',S}}} since
$H^0(X_{S'},\sO_{X_{S'}})=k$.
\end{proof}

\begin{prop}\label{prop2.8} All objects of $\sF(X_S)$ defined in Definition \ref{defn2.5} are vector
bundles. Thus $\sF(X_S)$ is an abelian, rigid $k$-linear tensor
category, the fiber functor 
\ga{2.17}{|_{x_S}: \sF(X_S)\to
{\rm Vect}_k, \ V\mapsto V|_{x_S}:=V\otimes_{\sO_{X_S}}k} defines
a Tannaka structure on $\sF(X_S)$.
\end{prop}
\begin{proof} For any $V\in\sF(X_S)$, we first show that there are $W_1,\,W_2\in\sC^N(X)$ and
a morphism $f:\pi_S^*W_1\to\pi_S^*W_2$ in $\text{Qcoh}(X_S)$ such
that $V=\text{coker}(f)$. One takes $W_2:=(\pi_S)_*V$ which by definition lies in $\sC^N(X)$, and defines $V_1$ to be the kernel of the surjection
$\pi_S^*W_2\surj V.$ Then $W_1:=(\pi_S)_*V_1\in\sC^N(X)$ and
$f:\pi_S^*W_1\surj V_1\inj \pi^*_SW_2$ satisfying
$\text{coker}(f)=V$.

Let $S'\subset\sC^N(X)$ be the full tensor subcategory generated
by $W_1$, $W_2$ and $S$. Then the pullbacks of $\pi_S^*W_1$ and
$\pi_S^*W_2$ to $X_{S'}$ (under $\pi_{S',S}:X_{S'}\to X_S$) become
trivial, thus $\pi_S^*W_1$ and $\pi_S^*W_2$ are in the image of
\ga{2.18}{\eta_{S',S}:\rm {Rep}(G_{S',S})\to \rm {Qcoh}(X_S)} By
Lemma \ref{lem2.6}, the functor $\eta_{S',S}$ is fully faithful,
thus $V=\text{coker}(f)$ is also in the image of $\eta_{S',S}$. In
particular, $V$ is a vector bundle.
\end{proof}

Next we show $\sF(X_S)=\sC^N(X_S)$ when $\pi(X,S,x)$ is reduced. To this aim, recall that a
bundle $V$ on $X$ is said to be \emph{strongly semistable} if for any nonsingular
projective curve $C$ and any  morphism $f:C\to X$, the
pullback $f^*V$ is semistable of degree 0 on $C$. It is to be noticed  that
 strongly semistable bundles on $X$ define a full tensor subcategory of ${\rm Qcoh}(X)$ which is $k$-linear, and that $\sC^N(X)$ is a full subcategory of it. 

\begin{lem}\label{lem2.9} If $\pi(X,S,x)$ is a smooth finite group scheme, and
if $V\in\sC^N(X_S)$, then $\pi_S^*(\pi_{S_*}V)\in\sC^N(X_S)$ and
$W:=\pi_{S_*}V$ is strongly semistable.
\end{lem}
\begin{proof} Let $G:=\pi(X,S,x)$, $\pi:=\pi_S$, $Y:=X_S$ and $y:=x_S$.
Consider 
\ga{2.19}{\xymatrix{\ar@/^20pt/[rr]^\phi Y\times_k G
\ar[r]^{\cong}\ar_{p_1}[dr] & Y\times_X Y\ar[r]\ar[d]
& Y\ar[d]^{\pi}\\
 & Y \ar[r]_{\pi}& X} }
where $\phi:Y\times_kG\to Y$ is the action of $G$, $p_1$ is the
projection to first factor and $Y\times_kG\cong Y\times_XY$ is
induced by $\psi=(p_1,\phi)$. Then $\pi^*\pi_*V\cong
{p_1}_*\phi^*V$. Let $k'\supset k$ be a finite Galois extension
such that $G_{k'}:=G\times_k\text{Spec}(k')=\sqcup_{g\in G(k^s)}
{\rm Spec} \ k'(g)$. Then 
\ga{2.20}{(\pi^*\pi_*V)\otimes
k'\in\sC^N(Y_{k'}).} 
As
$\pi^N(Y_{k'},y_{k'})\cong\pi^N(Y,y)\times_kk'$ (\cite[Chapter~II,
Proposition~5]{N2}), there exists $\sN\in\sC^N(Y)$ such that
$(\pi^*\pi_*V)\otimes k'$ is a subquotient of $\sN\otimes k'$. Let
$p:Y\times_k\text{Spec}(k')\to Y$ be the projection. Then
$(\pi^*\pi_*V)\otimes_{\sO_Y}p_*\sO_{Y_{k'}}$ is a subquotient of
$\sN\otimes_{\sO_Y}p_*\sO_{Y_{k'}}=\sN^{\oplus [k':k]}$, which
implies that
\ga{2.21}{(\pi^*\pi_*V)\otimes_{\sO_Y}p_*\sO_{Y_{k'}}=(\pi^*\pi_*V)^{\oplus
[k':k]}\in \sC^N(Y).} 
Thus $\pi^*\pi_*V\in\sC^N(Y)$.

To show that $W:=\pi_*V$ is strongly semistable, one has to consider 
curves $C$ mapping to $X$. This defines the fiber square
\ga{2.22}{\begin{CD}
Y_C @>g>> Y\\
@V\pi_C VV @V\pi VV\\
C @>f>> X
\end{CD}
}
So in particular, $\pi_C$ is still a principal bundle under 
$G$.
 One has
$f^*W={\pi_C}_*(V_C)$, where $V_C=g^*V$.
Since $V_C\in\sC^N(Y_C)$, we have just shown that
$\pi_C^*(f^*W)=\pi_C^*({\pi_C}_*V_C)\in\sC^N(Y_C)$. In particular,
$\pi_C^*(f^*W)$ is semistable of degree $0$, which implies that
$f^*W$ is semistable of degree $0$.
\end{proof}

\begin{thm}\label{thm2.10} When $\pi(X,S,x)$ is a smooth finite
group scheme, then $W:=\pi_{S_*}V\in\sC^N(X)$ whenever
$V\in\sC^N(X_S)$. In particular, $\sF(X_S)=\sC^N(X_S)$ and there
is an exact sequence of group schemes 
\ga{2.23}
{1\to\pi(X_S,x_S)\to\pi(X,x)\to\pi(X,S,x)\to 1}
\end{thm}

\begin{proof} By Lemma \ref{lem2.9}, $W:=\pi_{S_*}V$ strongly semistable and
$\pi_S^*W\in\sC^N(X_S)$. Let $\langle W\rangle \cup S$ (resp.
$\langle \pi_S^*W\rangle$) be the full tensor subcategory
generated by $W$ and objects of $S$ (resp. by $\pi_S^*W$).
Equipped with the fiber functor at $x$ (resp. at $x_S$), the
Tannaka category $\langle W\rangle \cup S$ (resp. $\langle
\pi_S^*W\rangle$) yields a Tannaka group scheme $G$ (resp.
$\pi(X_S,\langle \pi_S^*W\rangle,x_S)$). To show $W\in\sC^N(X)$,
it suffices to show that $G$ is a finite group scheme.

The full subcategory of $\langle W\rangle \cup S$ whose objects
become trivial when pulled-back to $X_S$ is precisely $S$ (see
\eqref{2.11}). The functor 
\ga{2.24}{\pi_S ^*:\langle W\rangle \cup
S\to \langle \pi_S ^*W\rangle} 
yields the sequence of
homomorphisms of group schemes 
\ga{2.25}{1\to
\pi(X_S,\langle\pi^*_SW\rangle,x_S)\to G\to \pi(X,S,x)\to 1 }
which we claim to be exact.

Indeed, the surjectivity of $G\to \pi(X,S,x)$ and the injectivity
of $\pi(X_S,\langle\pi_S^*W\rangle,x_S)\to G$ follow from the
definition and \cite[Proposition 2.21]{DeMil}. To show the
exactness at $G$, according to Theorem A.1, it is enough to check

\begin{itemize}\item[ (i)] An object $M$ of $\langle\pi_S ^*W\rangle$ is a quotient of an
object of the form $\pi_S^*N,\quad N\in \langle W\rangle \cup S$;
\item[(ii)] For an object $N$ of $\langle W\rangle \cup S$, the maximal
trivial subobject of $\pi_S^*N$ has the form $\pi_S^*N_0$ where
$N_0$ is a subobject of $N$.
\end{itemize}

The condition (i) is easy to see. By definition, $M$ is a
subquotient of $\pi_S^*N$, $N\in \langle W\rangle \cup S$. Thus
$\pi_{S_*}M$ is a subquotient of
$\pi_{S_*}\pi_S^*N=N\otimes\pi_{S_*} \pi^*_S\sO_{X_S}\in \langle
W\rangle \cup S$. Hence $\pi_{S_*}M$ lies in $\langle W\rangle \cup
S$. Now we have surjective map $\pi_S^*(\pi_{S_*}M)\to M$ and  we
can take $N=\pi_{S_*}M$.

As for (ii) we use projection formula
\ga{2.26}{H^0(X_S,\pi_S^*N)=H^0(X,\pi_{S_*}\pi_S^*N)\\ \notag
=\Hom_{\sO_X}(\pi_{S_*}\sO_{X_S}^\vee, N)= \bigoplus_{i=1}^rk\cdot
\phi_i} 
where $\phi_i:(\pi_{S_*}\sO_{X_S})^\vee\to N$. Let
$N_0=\sum_i\im(\phi_i)\subset N$. Then $N_0$ is in $S$ and any
morphism $\phi:(\pi_{S_*} \sO_{X_S})^\vee\to N$ has image in
$N_0$. By comparing the ranks,  we see that $\pi_S^*N_0$ is the
maximal trivial subbundle in $\pi_S^*N$.

Thus the sequence in \eqref{2.25} is exact. Since
$\pi(X_S,\langle\pi_S^*W\rangle,x_S)$ and $\pi(X,S,x)$ are finite
group schemes, so is $G$.

It is clear that Proposition \ref{prop2.8} implies
$\sF(X_S)\subset\sC^N(X_S)$ when $X_S$ is reduced. Then what we
have proved above implies $\sF(X_S)\supset\sC^N(X_S)$. The
sequence \eqref{2.23} follows from exactness in \eqref{2.25}
by taking projective limits on $\langle W\cup S\rangle$.
\end{proof}

The aim of our article is to understand the relationship between the groups
$\pi^N(X,x)$, $\pi^\et(X,x)$ and $\pi^F(X)$. Let us start with the following simple lemma.
\begin{lem} \label{lem2.11} The restriction homomorphisms
$\pi^N(X,x)\xrightarrow{r^{\acute{e}t}} \pi^{\acute{e}t}(X,x)$, resp.
$\pi^N(X,x)\xrightarrow{r^{F}} \pi^{F}(X,x)$
 induced by the inclusion of categories $\sC^{\acute{e}t}(X)\to \sC^N(X)$, resp. $\sC^F(X)\to \sC^N(X)$,
 with compatible fiber functor, are faithfully flat.
\end{lem}
\begin{proof}
By \cite[Proposition~2.21 (a)]{DeMil}, the homomorphism  $r$ (with
decoration) is faithfully flat if the inclusion of categories is
fully faithful, which is clear, and if any subobject in $\sC^N(X)$
of an object of the subcategory  is indeed a subobject in the
subcategory, which is also clear.
\end{proof}
\begin{nota} \label{nota2.12}
We set
$L(X,x)={\rm Ker} (\pi^N(X,x)\to \pi^{\acute{e}t}(X,x))$.
\end{nota}
\section{The representation category of the difference between Nori's  fundamental group scheme and it's \'etale part}
We continue to assume in this section that $X$ is a proper reduced
scheme defined over
 a field $k$ which is assumed to be perfect if of characteristic $p>0$. Most of the interest of this section is in the case when $k$ is a perfect field of characteristic $p>0$.  Further, 
$X$ 
is assumed to be reduced, connected in the sense that $H^0(X, \sO_X)=k$, and $X(k)\neq \emptyset$.
The purpose of this section is to determine the representation category
of $L(X,x)$ (see Notation \ref{nota2.12}).
To this aim, we first observe the following.
\begin{lem} \label{lem3.1}
 Let $S$ be a full tensor subcategory of $\sC^N(X)$ generated by finitely many objects.
Then
$$S^{\acute{e}t}:=S\cap \sC^{\acute{e}t}(X)$$
 is
 a full tensor subcategory generated by finitely many objects.
The fiber functor defined by $x$ identifies the Tannaka group of $S^{\acute{e}t}$
with the maximal reduced quotient of $\pi(X,S,x)$.
\end{lem}
\begin{proof}
Fullness is trivial, and again
 we apply \cite[Proposition~2.21 (a)]{DeMil}  to see that the Tannaka group
 $\pi((X,S^{\acute{e}t},x))$ of $S^{\acute{e}t}$ is a quotient of $\pi^N(X,S,x)$.
 In particular it is a finite group scheme over $k$, hence by Tannaka duality, $S^{\acute et}$ is
 finitely generated.

Tannaka duality shows that any quotient map $\pi(X,S,x)\surj H$ of
group schemes over $k$ yields a fully faithful functor from
$\text{Rep}(H)$ to $\sC^N(X)$ with image, say $(H)$, lying in $S$,
and consequently yielding an $H$-principal bundle $\pi_{(H)}:
X_{(H)}\to X$ which is proper, connected, with a rational
point $x_{(H)}$ mapping to $x$ so that $V\in (H) $ if and only if
$\pi_{(H)}^*(V)$ is trivial.  If $H$ is reduced then $(H)$
consists only of \'etale finite bundles, thus $(H)\subset S^\et$.
Hence $\pi(X,S,x)\to H$ factors through $\pi(X,S,x)\to
\pi(X,S^\et,x)\to H$. This shows that $\pi(X,S^\et,x)$ is the
maximal reduced quotient of $\pi(X,S,x)$.
\end{proof}

\begin{rmks} \label{rmks3.2}
1)
Any representation $V$ of an affine group scheme $G$ over $k$  can be embedded into
a direct sum of copies of the regular representation of $G$, i.e. the representation of $G$ on
its function algebra $k[G]$. Tannaka duality applied for a finitely generated tensor subcategory $S\subset \sC^N(X)$
shows that $S$ is generated by the sheaf $(\pi_S)_*(\sO_{X_S})$ on $X$.\\ \ \\
2)
In fact,
 $\pi^{\acute{e}t}(X,x)$ is defined by Grothendieck (\cite[Expos\'e~5]{Groth}) in the following sense.

Let us recall Grothendieck's definition. We assume as before that
$X$ is proper connected (i.e. $H^0(X, \sO_X)=k$) over $k$ perfect.
If $\bar{x}\in X$ is a geometric point, (that is $\bar{x}={\rm
Spec}(K)$ with $K$ algebraically closed), then Grothendieck
defines the fundamental group $\pi_1(X, \bar{x})$ as the
automorphism group of the fiber functor  $\sG r:=\{$finite \'etale
coverings of $ X$ defined over $k\} \to \{$Sets$\}, \ (\pi: Y\to
X)\mapsto \pi^{-1}(\bar{x})=\bar{x}\times_X Y$. He shows that it
is a pro-finite group, and  that $\sG r$ is equivalent to the
category of finite $\pi_1(X, \bar{x})$-finite sets. The natural
action of $\pi_1(X, \bar{x})$ on  a finite quotient $\pi_1(X,
\bar{x}) \surj H$ corresponds via the equivalence to a Galois
covering $\pi_H: X_H\to X$ with Galois group $H$. Recall here that
a Galois covering is a finite \'etale covering $\pi: Y\to X$ with
the property that $Y\times_XY$ is isomorphic over $k$ to a product
$Y\times_k H$, where $H$ is a constant group, and the constant
group $H$ acts on $Y$ with quotient $X$.

In order to compare the notions, let us first observe that if $k$
is perfect but not necessarly algebraically closed, and $X={\rm
Spec}(k)$ with the rational point point $x= X\in X(k)$, then $\sC^\et(X)$ is equivalent via
the fibre functor to ${\rm Vect}_k$ and consequently
$\pi^\et(X,x)=\{1\}$. On the other hand, Grothendieck's
fundamental group is then ${\rm Gal}(\bar{k}/k)$, which is highly
nontrivial. Thus we assume  $k=\bar{k}$, $X$ proper connected over
$k$,  $x=\bar{x}\in X(k)$ and $\bar{k}=K$. Then
\ga{}{\pi_1(X, x)=\varprojlim_{H  \ {\rm finite}, \ \pi_1(X, x)\surj H }H \notag \\
\pi^\et(X, x)=\varprojlim_{S \text{ with }  \pi(X,S,x)  \text{
finite reduced,} \atop \pi^N(X,x)\surj
\pi(X,S,x)}\pi(X,S,x).\notag} This means in particular that the
category $\sG r$ is uniquely determined by the Galois coverings
$\pi_H: X_H\to X$ and $\sC^\et(X)$ is uniquely determined by the
principal bundles $\pi_S: X_S\to X$ under a smooth finite group
over $\bar{k}=k$, and their transition morphisms. Thus given
$\pi(X,S,x)$ finite reduced, then $\pi(X,S,x)(k)$ is a  Galois
group. Vice-versa, given a Galois group $H$, understood as a
constant $k$-algebraic group, one has a surjection
$\pi^N(X,x)\surj H$ (\cite[I, Proposition 3]{N2}). This shows
\ga{}{\pi_1(X, x)=\pi^\et(X,x)(k). \notag } On the other hand, if
$k$ is perfect but not necessarily algebraically closed, $X$ is
proper connected over $k$ and $x\in X(k)$, it is quite elementary
to see that $\pi^\et(X,x)$ satisfies base change, i.e. the natural
$\bar{k}$-homomorphism \ga{}{\pi^\et(X\times_k \bar{k}, x\times_k
\bar{k}) \xrightarrow{\cong} \pi^\et(X,x)\times_k \bar{k}\notag}
is an isomorphism. Nori proves this with $\pi^\et(X,x)$ replaced by $\pi^N(X,x)$ in 
\cite[II, Proposition 5]{N2}, but his argument applies here (a fortiori) as well. 
 So we conclude in general 
\ga{}{ \pi_1(X \times_k
\bar{k}, \bar{x})=\pi^\et(X,x)(\bar{k}). \notag }

\end{rmks}

\begin{rmk} \label{rmk3.3}
The construction $S\mapsto [ \pi_{S}:(X_S, x_{S})\to (X,x)]$ is functorial in
the following sense. Given $S$ and  $T$
two finitely generated full tensor subcategories in $\sC^N(X)$, denote by $S\cup T$ the
category generated by objects of $S$ and $T$. Then it is also finitely
generated and Lemma \ref{lem2.6} shows that in the
following commutative diagram all maps are principal bundles
\ga{}{\begin{CD}(X_{S\cup T}, x_{S\cup T}) @>{\pi_{S\cup T, S}}>> (X_S, x_S)\\
@V{\pi_{S\cup T, T}}VV @VV\pi_SV\\
(X_T,x_T)@>>\pi_T> (X,x).
\end{CD}\notag}
Let now $S_1,S_2,S_3$ be finitely generated full tensor subcategories in $\sC^N(X)$.
Denote by $X_{ijk}$ the scheme $X_{S_i\cup S_j\cup S_k}$ and by $\pi_{ijk,ij}$ the map $\pi_{S_i\cup S_j\cup S_k,S_i\cup S_j}$.
Then we have the following commutative
cube
\ga{}{\xymatrix{
& X_{23} \ar@{->}'[d]^{\pi_{23,3}}[dd]
& &
X_{123}
\ar_{\pi_{123,23}}@{->}[ll]
\ar@{->}[dd]^{\pi_{123,12}}
\\
X_3
\ar^{\pi_{23,3}}@{<-}[ur]
\ar_{\pi_{13,3}}@{<-}[r]\ar@{-}[rr]
\ar_{\pi_3}@{->}[dd]
& & X_{13}
\ar^{\pi_{123,13}}@{<-}[ur]
\ar@{-}[d]^{\pi_{13,1}}\ar@{->}[dd]
\\
& X_2 \ar@{<-}'[r][rr]_{\pi_{12,2}}
& & X_{12}
\\
X
\ar@{<-}_{\pi_1}[rr]
\ar@{<-}^{\pi_2}[ur]
& & X_1 \ar@{<-}[ur]_{\pi_{12,1}}
}.\notag}
\end{rmk}

In the rest of this section, $S$ will denote a finitely generated
tensor full subcategory of $\sC^\et(X)$. Thus $\pi_S:X_S\to X$ is
\'etale and $X_S$ is reduced.

\begin{prop}\label{prop3.4}
Let $X$ be a proper reduced connected scheme defined over a
perfect field $k$ of characteristic $p>0$. Let $\pi_S: X_S \to X$
be an \'etale principal bundle, and let $V$ be an $F$-finite
bundle on $X$. Then if
$$\pi_S^*: H^0(X, \sO_X)\to H^0(X_S, \sO_{X_S})$$ is an isomorphism, so is
$$\pi_S^*: H^0(X, V)\to H^0(X_S, \pi_S^*V).$$
\end{prop}

\begin{proof} To simplify notation, let $Y:=X_S$, $\pi:=\pi_S$ in
the proof.

 Let $V_0$ be the maximal trivial subobject of $V$.
 Since $\pi$ is \'etale, the bundles associated to $\pi_*\sO_Y$, and therefore to
$(\pi_*\sO_Y)^\vee$,  are
 \'etale finite. The image under a morphism of $(\pi_*\sO_Y)^\vee$ to
$V$ is therefore at the same time \'etale- and $F$-finite, hence (see \eqref{2.12})
lies in $V_0$.
By projection formula we have
\ga{3.1}{H^0(Y,\pi^*V)=H^0(X, (\pi_*\sO_Y)\otimes V)\cong \Hom_{X}((\pi_*\sO_Y)^\vee,V)\\
\notag \subset \Hom_{X}((\pi_*\sO_Y)^\vee,V_0)\cong H^0(Y,\pi^*V_0).}
Hence
\ga{3.2}{H^0(X, V_0)\subset H^0(X, V)\subset H^0(Y, \pi^*V)\subset H^0(X, \pi^*V_0)\\
\notag
\stackrel{\rm (assumption)}= H^0(X, V_0),}
so one has everywhere equality.
\end{proof}

\begin{defn} \label{defn3.5} The category $\sD$ has for objects pairs $(X_S,V)$ where $S\subset \sC^\et(X)$ is a finitely generated full tensor subcategory,  $V\in\sC^F(X_S)$,
and for morphisms
$$\Hom((X_{S},V), (X_{S'}, W)):=\Hom_{X_{S\cup S'}}(\pi_{(S\cup S') , S }^*(V), \pi_{(S\cup S') ,S'}^*(W))$$
with notations as in Remark \ref{rmk3.3}.
\end{defn}

\begin{prop}  \label{prop3.6} The category
$\sD$ is an abelian, rigid $k$-linear tensor category with the following structures.
The composition rule  for
$$\varphi_{ij}\in {\rm Hom}((X_{S_i}, V_i), (X_{S_j}, V_j)),(i,j)=(1,2), (2,3)$$ is defined by
\ga{3.3}{\varphi_{13}=
\pi_{S_1\cup S_2\cup S_3, S_2\cup S_3 }^*(\varphi_{23})\circ  \pi_{S_1\cup S_2\cup S_3, S_1\cup S_2}^*(\varphi_{12})\notag \\
\in ({\rm Proposition} \ \ref{prop3.4}) \
{\rm Hom}((X_{S_1 }, V_1), (X_{S_3 }, V_3)).
}
 The additive structure is defined by
\ga{3.4}{
(X_{S}, V) \oplus (X_{T}, W)=(X_{S\cup T}, \pi_{S\cup T, S }^*(V)\oplus \pi_{S\cup T, T }^*(W))
}
and the tensor structure is defined by
\ga{3.5}{
(X_{S}, V) \otimes (X_{T}, W)=(X_{S\cup T}, \pi_{S\cup T, S }^*(V)\otimes \pi_{S\cup T, T }^*(W))
}
The functor
\ga{3.6}{\sD\xrightarrow{|_{x_{S }}} {\rm Vect}_k, \ (X_{S }, V)\mapsto V|_{x_{S }}} endows $\sD$ with the structure of a Tannaka category.

\end{prop}

\begin{proof} We check that $\sD$ is a category. The composition rule is well
defined as indicated thanks to Proposition \ref{prop3.4} which
implies
 \ga{3.7}{{\rm Hom}_{X_{S_1\cup S_2 \cup S_3}}
 ( \pi_{S_1\cup S_2\cup S_3 , S_1\cup S_2  }^*V_1,
\pi_{S_1\cup S_2\cup S_3 , S_2\cup S_3}^*V_3)= \\
{\rm Hom}_{X_{S_1\cup S_3}}(\pi_{S_1\cup S_3, S_1  }^*V_1,
\pi_{S_1\cup S_3,  S_3}^*V_3).\notag
}
Distributivity of the Hom with respect to the additive structure is obvious.
Then the unit object is
\ga{3.8}{\mathbb{I}=(X, \sO_X)}
 that is $S=\langle\sO_X\rangle$ and the endomorphism ring of the unit object is $k$.
 The Hom-sets are finite dimensional $k$-vector spaces. So $\sD$ is an additive category.
The dual object is defined by
\ga{3.9}{(X_{S }, V)^\vee=(X_{S }, V^\vee).}
 Since the $X_S$ is reduced, connected, and endowed with the canonical
 $k$-rational point, kernel, image and cokernel of a morphism in $\sD$ are still objects
 in $\sD$. This endows $\sD$ with an abelian structure, and then a tensor structure as well.
Furthermore, we observe
\ga{3.10}{(X_{S }, \sO_{X_{S }}) \mbox{ is  isomorphic  to }
(X, \sO_X) \mbox{ in }   \sD.}
More generally, for $S$ and $ T$ finitely generated full tensor subcategories of $\sC^\et(X)$, 
\ga{3.11}{(X_{S\cup T}, \pi_{S\cup T, S }^* V)
\mbox{ is isomorphic  to }  (X_{S }, V)\ \mbox{ in }   \sD  \ \rm{for \ all} \ T.
}
\end{proof}
\begin{lem} \label{lem3.7} Let $V\in\sC^N(X)$. Denote by $\langle V\rangle$ the
full tensor subcategory in $\sC^N(X)$ generated by $V$. Then the assignement
\ga{}{\sC^N(X)\xrightarrow{q} \sD, \ V\mapsto
(X_{\langle V\rangle^\et}, \pi_{\langle V\rangle^\et}^*(V))\notag}
defines a functor of k-linear rigid abelian tensor categories which commutes with the fiber functors $|_x$ and $|_{x_{\langle V \rangle^\et}}$.
\end{lem}
\begin{proof} The proof is obvious.
\end{proof}
We denote by $G(\sD)$ the Tannaka group scheme over $k$ defined by Proposition \ref{prop3.6}
via Tannaka duality. Then Lemma \ref{lem3.7} induces a homomorphism of group schemes
\ga{3.12}{q^*: G(\sD)\to \pi^N(X,x).}
\begin{thm} \label{thm3.8}
The homomorphism $q^*$ of \eqref{3.12} identifies $G(\sD)$ with $L(X,x)$
(from Notation \ref{nota2.12}). In other words, the representation category of $L(X,x)$ is equivalent to $\sD$.

\end{thm}
\begin{proof}
We follow the method of \cite[Section 5]{EP} which is recalled in the Appendix A. 
 The functor $q$ has the property that if $V$ is \'etale finite, then
$q(V)$ is isomorphic in $\sD$ to  $\oplus_1^r \mathbb{I}$ with
$r={\rm rank}(V)$. This implies that 
\ga{3.13}{G(\sD)\to L(X,s)}
and we have to show that \eqref{3.13} is an isomorphism. Let
$(X_{S},V)$ be an object in $\sD$. Then, by Theorem \ref{thm2.10},
$W:=(\pi_{S})_*V$
 is an object in $\sC^N(X)$. Moreover one has a surjection in $\sD$
\ga{3.14}{q(W)\surj (X_{S},V).} Thus every object of $\sD$ is a
quotient of the image by $q$ of an object in $\sC^N(X)$. Since we
can apply \eqref{3.14} to $(X_{S},V^\vee)$ as well, we conclude
that every object of $\sD$ is also a sub object of the image by
$q$ of an object in $\sC^N(X)$. This shows in particular that the
map \eqref{3.12} is injective, according to Theorem \ref{thmA1}. Thus the map in \eqref{3.13} is also injective.

We want to show that the map in \eqref{3.13} is also surjective. First notice
that if $V\in\sC^N(X)$ is such that $q(V)$ is trivial in $\sD$,  then $V$ is, in fact,
in $\sC^\et(X)$.
Furthermore, for $V\in\sC^N(X)$ let
 $A:\cong \oplus_1^r \mathbb{I}\inj q(V)$ be the maximal trivial subobject in $\sD$.
Then one has, applying projection formula and setting  $S:=\langle
V\rangle^\et$ 
\ga{3.15}{H^0(X_{S}, \pi_{S}^*(V))={\rm
Hom}_{\sO_X}(({\pi_S}_*\sO_{X_S})^\vee,
V)=\bigoplus_{i=1}^rk\cdot\varphi_i,} 
where
$\varphi_i:({\pi_S}_*\sO_{X_S})^\vee\to V$. Let $V_\et\subset V$
be the image of 
\ga{3.16}{\oplus^r_1\varphi_i:\quad
\bigoplus_1^r({\pi_S}_*\sO_{X_S})^\vee \to V.} As
$(\pi_{S*}\sO_{X_{S}})^\vee$ is \'etale finite, $V_\et$ is \'etale
finite and $\pi_S^*V_\et$ is a trivial bundle by \eqref{2.11}. On
the other hand, by the definition of $V_\et$, one has
\ga{3.17}{H^0(X, (\pi_{S*}\sO_{X_{S}})\otimes V)=H^0(X,
(\pi_{S*}\sO_{X_{S}})\otimes V_\et).} Applying the projection
formula again, one has 
\ga{3.18}{H^0(X,
(\pi_{S*}\sO_{X_{S}})\otimes V_\et)=H^0(X_{S}, \pi_{S}^*(V_\et)).}
Thus we conclude that 
\ga{3.19}{A= (X_S,\pi_{S}^*(V_\et)).} This
shows that the maximal trivial subobject of $q(V)$ is equal to
$q(V_\et)$.

Now we are in the situation of the assumption of
Theorem \ref{thmA1},
and we conclude that \eqref{3.13} is an isomorphism. \end{proof}

\section{Nori's fundamental group scheme surjects onto the product of its \'etale part with its local part}
The aim of this section is to prove the following
\begin{thm} \label{thm4.1} Let $X$ be a reduced proper scheme defined over a perfect field $k$ of
characteristic $p>0$, which is connected, i.e. $H^0(X, \sO_X)=k$, and let $x\in X(k)$. Then the group scheme homomorphism
$$\pi^N(X,x)\xrightarrow{r^\et \times r^F} \pi^\et(X,x)\times \pi^F(X,x)$$ over $k$ defined in Lemma \ref{lem2.11} is faithfully flat.
Equivalently, the group scheme homomorphism $$L(X,x)\xrightarrow{r^F|_{L(X,x)}}
\pi^F(X,s)$$ over $k$ is faithfully flat.
\end{thm}

We will need the following theorem, which is a generalization of
Proposition \ref{prop3.4} where $W=\sO_{X_S}$.
\begin{thm} \label{thm4.2}
Let $X$ be a proper reduced connected scheme defined over a
perfect field $k$ of characteristic $p>0$. Let $\pi_S: X_S\to X$
be a principal bundle under an \'etale finite group scheme $G$,
let $V$ be a
 $F$-finite bundle over $X$, and let $\varphi: W\inj \pi^*V$
be an inclusion of an $F$-finite bundle $W$ on $X_S$. Then if
$\pi_S^*: H^0(X, \sO_X)\to H^0(X_S, \sO_{X_S})$ is an isomorphism,
there exists an $F$-finite bundle $V_0$ on $X$, together with an
inclusion $i: V_0\inj V$ so that $\varphi=\pi^*i: \pi^*V_0=W\inj
\pi^*V$.
\end{thm}
\begin{proof} To simplify notation, we set $Y:=X_S$, $\pi:=\pi_S$ in the proof.

We first assume that $W$ is simple. The bundle $V$ has a
decomposition series $V_0\subset V_1\subset \ldots \subset V_N=V$.
This means that the graded pieces $V_i/V_{i-1}$ are simple. Assume
$W\neq 0$, so $\varphi \neq 0$. Let $M\le N$ be the largest index
for which the image of $\varphi(W)$ in $\pi^*(V_M/V_{M-1})$ is not
trivial. Since $W$ is simple, then $\varphi$ induces an injective
morphism $W\inj \pi^*(V_M/V_{M-1})$. Since $V_M/V_{M-1}$ is
$F$-finite as well, we may assume that $V$ itself is simple. Then
we have to show that
 if $\varphi\neq 0$, then it is an isomorphism and we take $i={\rm identity}$.
Further we can replace $W$ by its maximal subobject which is simple. So both $W$ and $V$ are simple now.
Projection formula says
\ga{4.1}{H^0(Y, W^\vee \otimes \pi^*V)=H^0(X, \pi_*(W^\vee)\otimes V) \ {\rm equivalently}\\
{\rm Hom}(W, \pi^*V)={\rm Hom}((\pi_*(W^\vee))^\vee, V).\notag}
So
\ga{4.2}{\varphi \in {\rm Hom}(W, \pi^*V) \leftrightsquigarrow
\psi \in {\rm Hom}((\pi_*(W^\vee))^\vee, V).}
The projection $\pi^*(\pi_*W^\vee) \surj W^\vee$ is dual to the injection
\ga{4.3}{\iota: W\inj (\pi^*(\pi_*W^\vee))^\vee}
and one has commutative diagram
\ga{4.4}{\xymatrix{
(\pi^*(\pi_*W^\vee))^\vee \ar^{\quad\pi^*\psi}[rr] && \pi^*V\\
\ar_\iota[u] \ar_\phi[rru]
W} 
}
On the other hand, as $\pi$ is principal bundle under a finite smooth group scheme, one has the fiber square diagram
\ga{4.5}{
\xymatrix{\ar@/^20pt/[rr]^\mu
Y\times_k G \ar[r]_{\theta \cong}\ar[d]_{p_1}
& Y\times_X Y \ar[r]_{p_2}\ar[d]_{p_1}
& Y\ar[d]^\pi\\
 Y \ar@{=}[r]& Y \ar[r]_\pi& X
} } where $\mu$ is the action of $G$ on $Y$. Flat base change
theorem implies
\ga{4.6}{(\pi^*(\pi_*W^\vee))=(p_1)_*(\mu^*W^\vee).} We now make
base change from $k$ to its separable closure $\bar{k}$ (which is its
algebraic closure as $k$ is perfect). Then 
\ga{4.7}{G\times_k
\bar{k}=\sqcup_{g\in G(\bar{k})} {\rm Spec} \ \bar{k}(g)} and consequently
\ga{4.8}{(p_1)_*(\mu^*W^\vee)\times_k \bar{k}=\oplus_{g\in G(\bar{k})} W^\vee_g \ {\rm with}\\
W^\vee_g=\mu^*(W^\vee)|_{Y\times_k {\rm Spec} \ \bar{k}(g)}.\notag}
Now we assume $\varphi\neq 0$. Since $V$ is simple, and since by
Theorem \ref{thm2.10}, 
$(\pi_*W^\vee)^\vee$ is Nori finite on $X$, we conclude that
$\psi$ 
\ga{4.9}{\psi: (\pi_*W^\vee)^\vee \surj V} is surjective.
This implies that 
\ga{4.10}{\pi^*\psi:
\pi^*(\pi_*W^\vee)^\vee\surj \pi^*V} is surjective as well. So we
conclude that over $k^s$ one has
\ga{4.11}{\pi^*\psi\times_k \bar{k}: \oplus_{g\in G(\bar{k})} W_g \surj \pi^*V\times_k \bar{k}\\
W_1=W,  \ \pi^*\psi\times_k \bar{k}|_W=\varphi.\notag}
On the other hand,  since the Nori's fundamental group $\pi^N(X,x)$
respects base change by a separable extension of $k$ (see \cite[II, Proposition 5]{N2}),
$W$ being simple over $Y/k$ implies that $W\times_k \bar{k}$
remains simple.
We conclude that in $\sC^F(X\times_k \bar{k})$
the object $\pi^*V\times_k \bar{k}$ is a quotient of a sum of simple objects, thus is a sum of simple objects itself. Thus we have
\ga{4.12}{W\oplus
\bigoplus_{{\rm some}\  g\in G(\bar{k})\setminus 1} W_g\xrightarrow{ \cong} \pi^*V\times_k \bar{k}.}
On the other hand, one has
\ga{4.13}{{\rm End}(\pi^*V\times_k \bar{k})={\rm End}(\pi^*V)\otimes_k \bar{k}\\
\stackrel{(\text{Prop.} \ref{prop3.4})}={\rm End}(V)\otimes_k \bar{k}
\stackrel{(\text{as $V$ simple})}= k\otimes_k \bar{k}=\bar{k}.\notag}
Consequently 
\ga{4.14}{\bigoplus_{{\rm some}\  g\in
G(\bar{k})\setminus 1} W_g=0} and we conclude
\ga{4.15}{\varphi\times_k \bar{k}={\rm isomorphism}, \ {\rm thus} \
\varphi \ {\rm isomorphism}.} Summarizing, we see
\ga{4.16}{\varphi: W\to \pi^*V \ {\rm as \ in \ Theorem \ with} \ W \ {\rm simple}  \ \Longrightarrow\\
\exists V_0\subset V, \ V_0 \ {\rm simple, \ with}\notag\\
\big(\varphi(W)\subset \pi^*V\big)=\big(\pi^*(V_0)\subset \pi^*(V)\big).\notag}
We treat now the general case. Let $W_0\subset W$ be a proper subbundle of
maximal rank so that $W_1=W/W_0$ is simple. Since $\varphi|_{W_0}: W_0\inj \pi^*V$,
induction on the length of a decomposition series for $W$ says there is
an $F$-finite subobject $U_0\subset V$ with $W_0=\pi^*(U_0)$ and the inclusion
$\varphi(W_0)\subset V$ being the inclusion $\pi^*(U_0)\subset \pi^*(V)$.
This implies that the induced map $ \varphi_1:W_1 \to \pi^*(V/U_0)$ is injective as well.
By \eqref{4.16}, we know that there is an $F$-finite bundle $U_1\subset V/U_0$ so that
the inclusion $\varphi_1(W_1)\subset \pi^*(V/U_0)$ is the inclusion $\pi^*(U_1)\subset
\pi^*(V/U_0)$. Denoting by $\tau: V\to V/U_0$ the projection, we conclude that
$V_0=\tau^{-1}(U_1)$ has the property that the inclusion $\pi^*(V_0)\subset \pi^*(V)$
 is precisely the inclusion
$\varphi(W)\subset \pi^*V$. This finishes the proof.
\end{proof}

\begin{proof}[Proof of Theorem \ref{thm4.1}]
Notice that the functor $q|_{\sC^F(X)}:{\sC^F(X)}\to \sD$, which is Tannaka dual
to $r^F|_{L(X,x)}$, is particularly easy to understand. To an $F$-finite
bundle $V$ on $X$, it assigns the pair $(X,V)$ as $X=X_{\langle \sO_X \rangle}$.
By \cite[Proposition~2.21 (a)]{DeMil}, the homomorphism
$r^F|_{L(X,s)}: L(X,x)\to \pi^F(X,x) $ is faithfully flat if the inclusion of
categories is  fully faithful, which is
then clear, and if any subobject in $\sD$ of an object coming from  the
subcategory $\sC^F(X)$  is indeed coming from the subcategory. This
is precisely Theorem \ref{thm4.2}.
\end{proof}
\begin{rmk}\label{rmk4.3}
By \cite{R}, Th\'eor\`eme 4.3.1, if $X$ is a smooth projective curve of genus
$g\ge 2$ over an algebraic closed  field $k$ of characteristic $p>0$, then for
$\ell \neq p $ prime with
$\ell+ 1 \ge (p-1)g$, there is a cyclic covering $\pi: Y\to X$ of degree $\ell$
(thus \'etale) so that ${\rm Pic}^0(Y)/{\rm Pic}^0(X)$ is ordinary.
Since $\pi$ is Galois cyclic of order $\ell$, it is defined as
${\rm Spec}_X(\oplus_0^{\ell-1} L^i)$ for some $L$ \'etale finite of
rank 1 over $X$ and of order $\ell$, thus $\pi=\pi_{\langle L\rangle}$ and $Y=X_{\langle L\rangle}$.
We conclude that there are $p$-power torsion rank 1 bundles on $X_{\langle L\rangle}$
which do not come from $X$. Consequently,
$${\rm Ker}(r^\et\times r^F)={\rm Ker} (r^F|_{L(X,x)})$$ is not trivial in general.
The purpose of the next section is to determine the structure of its representation category.
\end{rmk}
\section{The quotient category}
The aim of this section is to describe the representation category
of the kernel of the homomorphism  $r^\et\times r^F$ studied  in Theorem \ref{thm4.1}.
Recall that $X$ is a reduced proper scheme over a perfect field of
characteristic $p>0$ with a  rational point $x\in X(k)$
and is connected in the sense that  $H^0(X, \sO_X)=k$.

We first need some notations. Consider the inductive system
$\sT^\ell$ of finitely generated tensor  full subcategories
$T\subset \sC^F(X)$, which yields the projective system  of
$k$-schemes 
\ga{5.1}{X_{\sT^\ell}:=\varprojlim_{T\in \sT^\ell}
X_T} with transition maps for $T\subset T'$
\ga{5.2}{\pi_{T', T}:
X_{T'}\to X_T.} Similarly we define  the inductive system
$\sS^\et$ of finitely generated full tensor subcategories
$S\subset \sC^\et(X)$.

Let $S\subset \sC^N(X)$ be a finitely generated tensor full
subcategory and let $\pi_S: X_S\to X$ be the corresponding principal
bundle as in \eqref{2.7}. For any bundle $V$ on $X_S$, we define
the inductive system of finite dimensional $k$-vector spaces
\ga{5.3}{H^0_{\sT^\ell}(X_S, V)=\varinjlim_{T\in \sT^\ell}
H^0(X_{S\cup T}, \pi^*_{S\cup T, S}(V)),} where the structure map
for an inclusion $T\subset T'$ is defined by 
\ga{5.4}{H^0(X_{S\cup
T}, \pi^*_{S\cup T, S}(V)) \xrightarrow{\pi^*_{S\cup T', S\cup T}}
H^0(X_{S\cup  T'}, \pi^*_{S\cup  T', S}(V)).}

\begin{lem} \label{lem5.1} For any $V\in\sF(X_S)$ (see Definition \ref{defn2.5}),
the $k$-vector space $H^0_{\sT^\ell}(X_S, V)$ is finite
dimensional, and is endowed with a $k$-linear embedding $
H^0_{\sT^\ell}(X_S, V) \subset V|_{x_S}$.\end{lem}
\begin{proof}
Recall from  \eqref{2.10} that one has the constant projective
system of rational points $x_{S\cup T}\in X_{S\cup T}(k)$. Using
Proposition \ref{prop2.8}, it is easy to see $\pi^*_{S\cup
T,S}(V)\in\sF(X_{S\cup T})$ and 
\ga{5.5}{H^0(X_{S\cup T},
\pi^*_{S\cup T, S} (V))\subset \pi^*_{S\cup T, S} (V)|_{x_{S\cup
T}}=V|_{x_S}.} 
Since the structure map in \eqref{5.4} of our
inductive system are all injective, the claim of Lemma
\ref{lem5.1} follows from \eqref{5.5}.
\end{proof}

In the rest of the section, $S$ will denote a finitely generated
tensor full subcategory of $\sC^\et(X)$ and $T$ will denote a
finitely generated tensor full subcategory of $\sC^F(X)$.

\begin{defn}\label{defn5.2} For $S$ 
a finitely generated
tensor full subcategory of $\sC^\et(X)$ and $T$ a
finitely generated tensor full subcategory of $\sC^F(X)$, one defines
$\sE(X_{S\cup T})\subset\sF(X_{S\cup T})$
to be the full subcategory whose objects $V$ have the property
that $(\pi_{S\cup T,S})_*V\in\sC^F(X_S)$.
\end{defn}

\begin{lem}\label{lem5.3}
Let $S\subset S'$ be in  $\sS^\et$ and $T\subset T'$ be
in $\sT^\ell$, and $V$ be an object in $\sE(X_{S\cup
T})$. Then
\begin{itemize}
\item[1)] The following commutative diagram is cartesian \ga{}{\begin{CD}
X_{S'\cup T} @>\pi_{S'\cup T, S'}>> X_{S'}\\
@V\pi_{S'\cup T,S\cup T}VV @VV\pi_{S',S} V \\
X_{S\cup T} @>>\pi_{S\cup T,S} > X_S
\end{CD}\notag}
\item[2)] $\sE(X_{S\cup T})$ is a k-linear abelian, rigid tensor
category.
\item[3)]
 $\pi^*_{S\cup T',S\cup T}V\in\sE(X_{S\cup T'})$.
\item[4)]
 $\pi^*_{S'\cup T,S\cup T}V\in\sE(X_{S\cup T})$.
\item[5)]
 the canonical homomorphism 
\ga{5.6}{H^0_{\sT^\ell}(X_{S\cup T},
V)\to H^0_{\sT^\ell}(X_{S'\cup T}, \pi^*_{S'\cup T, S\cup T}(V))}
is an isomorphism. \end{itemize}
\end{lem}
\begin{proof} 1) We first show the following commutative diagram 
\ga{5.7}{\begin{CD}
X_{S\cup T} @>\pi_{S\cup T, S}>> X_S\\
@V\pi_{S\cup T, T}VV @VV\pi_S V \\
X_T @>>\pi_T > X
\end{CD}}
is  cartesian. Indeed, $\pi_{S}:X_S \to X$ is a principal bundle
under $\pi(X, S, x)$, and similarly for $\pi_T, \pi_{S\cup T}$. So
the assertion is equivalent to show that the natural homomorphism
\ga{5.8}{\pi(X, S\cup T, x)\to \pi(X,S,x)\times \pi(X, T, x)}
induced by the embeddings $S\subset (S\cup T), T\subset (S\cup T)$
of categories is an isomorphism. By Theorem \ref{thm4.1} and the
profinite property of the various Tannaka group schemes involved,
the homomorphism
$$\pi^N(X,x)\to \pi(X,S,x)\times \pi(X, T, x)$$ is
faithfully flat. Thus  so is its factorization  \eqref{5.8}. On
the other hand, by definition, every object in $S\cup T$ is a
subquotient of tensors of objects in $S$ and objects in $T$, thus
by \cite[Proposition 2.21]{DeMil}, \eqref{5.8} is injective.
Therefore we have 
\ga{5.9}{X_{S'}\times_{X_S}X_{S\cup
T}=X_{S'}\times_{X_S}(X_S\times_X X_T)=X_{S'\cup T}.} That is we
have the cartesian product 
\ga{}{\begin{CD}
X_{S'\cup T} @>\pi_{S'\cup T, S'}>> X_{S'}\\
@V\pi_{S'\cup T,S\cup T}VV @VV\pi_{S',S} V \\
X_{S\cup T} @>>\pi_{S\cup T,S} > X_S
\end{CD}\notag}
of Lemma \ref{lem5.3}, 1).

2) Note that $\pi_{S\cup T,S}:X_{S\cup T}\to X_S$ is a principal
$\pi(X,T,x)$-bundle since \eqref{5.7} is cartesian, thus
$(\pi_{S\cup T,S})_*\sO_{X_{S\cup T}}$ is $F$-finite on $X_S$.  The
pullback by $\pi_{S\cup T,S}$ of any $F$-finite bundle on $X_S$
lies in $\sE(X_{S\cup T})$. Then it is easy to write any
$V\in\sE(X_{S\cup T})$ as a cokernel of a morphism $\pi_{S\cup
T,S}^*W_1\to \pi_{S\cup T,S}^*W_2$, where
$W_1,\,W_2\in\sC^F(X_S)$, thus  2) follows.

3) and 4) are easy by playing diagram game and using projection
formula, as we did already many times. 

To show 5), one use Proposition \ref{prop3.4} and projection
formula 
\ga{5.10}{H^0(X_{S'\cup T}, \pi^*_{S'\cup T,S\cup T}V)=
H^0(X_{S'}, (\pi_{S'\cup T,S'})_*\pi^*_{S'\cup T,S\cup T}V)\\
=H^0(X_{S'}, \pi^*_{S',S}(\pi_{S\cup T,S})_*V)\stackrel{(\text{Prop.
} \ref{prop3.4})}=
H^0(X_S,(\pi_{S\cup T,S})_*V)\notag\\
\notag = H^0(X_{S\cup T}, V)}
Of course \eqref{5.10} remains valid
under base change $X_{S\cup T'}\to X_{S\cup T}$ for any
$T'\in\sT^\ell$ containing $T$. The claim follows when taking
limit on $T'$.
\end{proof}

\begin{defn}\label{defn5.4}
We define the category $\sE$ by its objects and its morphisms.
\begin{itemize}
\item[1)] The objects of $\sE$ are pairs $(X_{S\cup T}, V)$
where $S\in \sS^\et$, $T\in \sT^\ell$ and $V\in\sE(X_{S\cup T})$.
\item[2)] The morphisms of $\sE$ are defined as follows. Set $S'':=S\cup S'$ and
$T'':=T\cup T'$ and define
\ga{}{{\rm Hom}((X_{S\cup T}, V), (X_{S'\cup T'}, V')):=\notag \\
  H^0_{\sT^\ell}(X_{S''\cup T''}, \pi_{S'' \cup T'', S\cup T}^*(V^\vee)
 \otimes  \pi_{S''\cup T'', S'\cup T'}^*(V') ).\notag}
\end{itemize}
\end{defn}

\begin{prop} \label{prop5.5} The category
$\sE$ defined in Definition \ref{defn5.4} is a $k$-linear abelian,
rigid tensor category, with tensor product defined by
\ga{5.11}{ (X_{S\cup T}, V) \otimes  (X_{S'\cup T'}, V'):=\\
(X_{S'' \cup T''},
\pi_{S'' \cup T'', S\cup T}^*(V)\otimes \pi_{S'' \cup T'', S'\cup T'}^*(V'))\notag}
($S'':=S\cup S'$, $T'':=T\cup T'$) and $\mathbb{I}$-object
\ga{5.12}{\mathbb{I}:=(X, \sO_X).}
It has a fiber functor
\ga{5.13}{\sE\to {\rm Vect}_k, \ \ (X_{S\cup T}, V)\mapsto V|_{x_{S\cup T}}}
defining a Tannaka group scheme $G(\sE,x)$ over $k$.
\end{prop}
\begin{proof}
We first show that $\sE$ is a category.
We just have to see that we can compose morphisms. Let
$\varphi: A\to A'$ be a morphism in $\sE$.
 Then according to Lemma \ref{lem5.3} we can write
 $A=(X_{S\cup T}, V)$ and $A'=(X_{S\cup T'},V')$ and find
a category $T''\in \sT^\ell$, containing both $T$ and $T'$,  such that $\varphi$
is represented by a morphism
 \ga{5.14}{\varphi_{T''}\in \Hom(\pi^*_{S \cup T'', S\cup T}(V),
\pi^*_{S \cup T'', S\cup T'}V')).}
It is easy to see that this presentation allows us to define the
composition of morphisms in $\sE$.  The kernel and image of $\varphi_{T''}$
are defined to be $X_{S\cup S'\cup T''}$ equipped with the
kernel and image of $\varphi_{T''}$, respectively. One might notice that
the kernel (image) of $\varphi$ does depend on the choice of $T''$.
However notice that for an object $(X_{S\cup T},V)$ in $\sE$, and for any
$T'\in \sT^\ell$, $T'\supset T$
\ga{5.15}{(X_{S\cup T},V) \ {\rm is \ isomorphic \ with} \
(X_{S\cup T'}, \pi_{S\cup T',S\cup T}^*(V)) \ \text{ in } \ \sE.}
It is clear how to define the tensor product in $\sE$.
The dual of an object is defined by
\ga{5.16}{(X_{S\cup T}, V)^\vee=(X_{S\cup T}, V^\vee).}
\end{proof}
\begin{defn} \label{defn5.6}
We define the functor
$q: \sD\to \sE, (X_S, V)\mapsto (X_S, V)$ which is clearly compatible with the fiber
 functors of $\sD$ and $\sE$.
\end{defn}
\begin{thm} \label{thm5.7}
The $k$-group scheme homomorphism $$q^*: G(\sE,x)\to L(X,x)$$
induced by $q$ is isomorphic to  $\text{Ker}(L(X,x)\to
\pi^F(X,x))$ and consequently is isomorphic to
$${\rm Ker}
\left(r^{\et}\times r^F:\pi^N(X,x)\to \pi^\et(X,x)\times
\pi^F(X,x)\right).$$
 In other words the representation category of
${\rm Ker}(r^\et\times r^F)$ is equivalent to $\sE$.
\end{thm}
\begin{proof}
We follow again the method in \cite[Section 5]{EP} which is recalled in the Appendix A. By means of
Theorem \ref{thm4.2}, $\sC^F(X)$ is a full subcategory of $\sD$.
Furthermore, if $V\in\sC^F(X)$, then its pull back on $X_{\langle
V\rangle}$ is trivializable, thus via the composite functor
\ga{5.17}{\sC^F(X)\to \sD\xrightarrow{q} \sE } every object in
$\sC^F(X)$ is trivializable. This implies the factorization
\ga{5.18}{G(\sE)\to {\rm Ker} \left(L(X,x)\to \pi^F(X,x)\right)}
of $q^*$. We have to show that \eqref{5.18}  is an isomorphism.

If $(X_{S\cup T}, V)$ is an object in $\sE$, then
\ga{5.19}{\pi_{S\cup T, S}^*(\pi_{S\cup T, S})_* V\surj V} and
since $(X_S, (\pi_{S\cup T, S})_* V)$ is an object of $\sD$, every
object of $\sE$ is the quotient of an object coming from  $\sD$
via $q$. Thus, according to Theorem \ref{thmA1}, \eqref{5.18}  is injective. In order to conclude  that \eqref{5.18}
is an isomorphism, one is led to show the
following facts:
\begin{itemize}\item[ (i)] An object $A=(X_S,V)$ in $\sD$ for which $q(A)$ becomes
trivial in $\sE$ comes from $\sC^F(X)$, that is $V\cong \pi_S^*(W)$ for
some $W\in \sC^F(X)$
\item[(ii)] Let $B$ be the maximal trivial subobject of $q(A)$ in $\sE$. Then there
exists a subobject $A_0$ of $A$ in $\sD$ such that $B\cong q(A_0)$ in $\sE$.
\end{itemize}
In fact, by definition (see the proof of Proposition \ref{prop5.5}),
$B$ is defined in terms of a category $T\in \sT^\ell$ as the
pair $(X_{S\cup T},V_0)$ where $V_0$ is the maximal trivial subobject of $\pi_{S\cup T,S}^*(V)$
in $\sE(X_{S\cup T})$.  Proposition \ref{prop5.8} below (which holds in fact for any $T$),
shows that there exists an $F$-finite bundle $W\in T$ on $X$ and an inclusion $j:\pi_{S}^*W\inj V$
such that $$\pi_{S\cup T,S}^*(j):\pi_{S\cup T}^*W\cong V_0.$$
This implies both properties (i) and (ii) above.
\end{proof}

\begin{prop} \label{prop5.8}
Let $V\in\sC^F(X_S)$, and $V_0\subset \pi^*_{S\cup T,S} (V)$ be
the maximal  trivial subobject of $\pi^*_{S\cup T,S} (V)$ in
$\sE(X_{S\cup T})$. Then there exists an $F$-finite bundle $W\in
T$ on $X$ equipped with an inclusion $j: \pi^*_S(W)\inj V$ such
that
$$\pi^*_{S\cup T,S}(j): \pi_{S\cup  T}^*(W)\cong V_0\subset \pi_{S\cup T, S}^*(V).$$
\end{prop}
\begin{proof}
Using  the cartesian diagram \eqref{5.7}, we have
\ga {5.20}{H^0(X_{S\cup T}, \pi_{S\cup T,S}^*V)\cong
H^0(X_T,\pi_{S\cup T,T*}\pi_{S\cup T,S}^*V) \\\notag
\stackrel{\eqref{5.7}}\cong H^0(X_T,\pi_{T}^*\pi_{S*}V)\cong
H^0(X,{\pi_T}_*\sO_{X_T}\otimes{\pi_S}_*V)\\ \notag
\cong \Hom_X(({\pi_T}_*\sO_{X_T})^{\vee}, {\pi_S}_*V)
=\oplus_1^rk\cdot\varphi_i
}
 for some morphisms $\varphi_i: ({\pi_T}_*\sO_{X_T})^\vee \to{\pi_S}_* V$.
 Let $W\subset {\pi_S}_* V$ be the image of
\ga{5.21}{ \oplus_1^r\varphi_i: \oplus_1^r
({\pi_T}_*\sO_{X_T})^\vee \to {\pi_S}_*V.} Then $W$ is
trivializable by $\pi_T$ and in particular $W\in T$ is $F$-finite.
We show that the map $j:\pi^*_SW\to V$ induced from the inclusion
$i:W\inj \pi_{S*}V$ is injective. Indeed, by definition,
$\pi_{S*}j$ composed with the canonical injection $W\to
\pi_{S*}\pi_S^*W$
gives back $i$ 
\ga{5.22}{\xymatrix{W\ar^{}[d] \ar^{i}[dr]\\
\pi_{S*}\pi^*_SW\ar_{\pi_{S*}j}[r] & \pi_{S*}V }} Let $V'$ be the
image of $j$ in $V$. Thus $V'$ is a quotient of the pull-back by
$\pi_S$ of the $F$-finite bundle $W$ on $X$ and according to
Theorem \ref{thm4.2} applied to $W^\vee$, we conclude that there
is a quotient $q:W\to W'$ of $W$ such that $V'\cong \pi_S^*W'$. By
the diagram in \eqref{5.22} applied to $W'$ there exists map
$i':W'\to \pi_{S*}V$ inducing the inclusion $\pi^*_SW'\cong V'\inj
V$. Thus the diagram in \eqref{5.22} factors as follows
\ga{5.23}{\xymatrix{ W\ar[r]^q\ar[d]\ar@/^40pt/[rrd]^i&W'\ar^{}[d] \ar^{i'}[dr]\\
\pi_{S*}\pi^*_SW\ar[r] & \pi_{S*}\pi^*_SW'\ar_{\pi_{S*}j}[r]
& \pi_{S*}V
}}
Since $q$ is epimorphism and $i$ is monomorphism, we conclude that $q$ is in fact an
isomorphism. Thus $j: \pi_S^*W\to V$ is in fact injective.

On the other hand, one has from \eqref{5.21}
\ga{5.24}{\Hom_X(({\pi_T}_*\sO_{X_T})^{\vee},
{\pi_S}_*V)=\Hom_X(({\pi_T}_*\sO_{X_T})^{\vee}, W)}
 Thus
 \ga{5.25}{H^0(X_{S\cup T},
\pi_{S\cup T,S}^*V)=H^0(X_T,\pi_T^*W)=H^0(X_{S\cup T},\pi_{S\cup T}^*W)}
 which means that $\pi^*_{S\cup T}(W)=V_0\subset\pi^*_{S\cup T,S}(V)$.
\end{proof}
\begin{rmk} \label{rmk5.9}
Using Proposition \ref{prop2.8}, replacing $X$ by $X_S$ and $X_{S}$ by 
$X_{S\cup T}$, one sees that the objects of $\sE(X_{S\cup T})$ are precisely those bundles which come from a representation of a local 
fundamental group over $k$. 
\end{rmk}
\begin{appendix} \section{Exact sequence of Tannaka group schemes}
In this appendix, we summarize the material on Tannaka categories which we used throughout the article. The statements and their proofs are taken from  \cite{EP}, but for the reader's convenience, we gather the information in a compact form here.

Let $L\xrightarrow q G\xrightarrow p A$ be a sequence of
homomorphism of affine group scheme over a field $k$. It induces a
sequence of functors
\ga{a1}{\Rep(A)\xrightarrow{p^*}\Rep(G)\xrightarrow{q^*}\Rep(L)}
where $\Rep$ denotes the category of finite dimensional
representations over $k$.
\begin{thm} \label{thmA1} With the above settings we have
\begin{itemize}
\item[(i)] the map $p:G\to A$ is faithfully flat (and in particular
surjective) if and only if $p^*$ is a full subcategory of $\Rep(G)$, closed
under taking subquotients.
\item[(ii)] the map $q:L\to G$ is a closed immersion if and only if any object of
$\Rep(L)$ is a subquotient of object of the form $q^*(V)$ for some $ V\in
\Rep(G)$.
\item[(iii)] Assume that $q$ is a closed immersion and that $p$ is
faithfully flat. Then the sequence $L\xrightarrow q G\xrightarrow
p A$ is exact if and only if the following conditions are fulfilled: \end{itemize}
\begin{itemize}\item[(a)] For an object $V\in \Rep(G)$, $q^*(V)$ in $\Rep(L)$ is trivial if and only if $V\cong p^*U$
for some $U\in \Rep(A)$.
\item[(b)] Let $W_0$ be the maximal trivial subobject of $q^*(V)$ in $\Rep(L)$. Then there exists
$V_0\subset V$ in $\Rep(G)$, such that $q^*(V_0)\cong W_0$.
\item[(c)] Any $W$ in $\Rep(L)$ is embeddable (hence by taking duals a quotient) of $q^*(V)$ for
some $V\in \Rep(G)$.
\end{itemize}
\end{thm}
\begin{proof}
The statements (i) and  (ii) are due to Deligne-Milne. We refer to \cite[Proposition
2.21]{DeMil}. We show (iii).

Assume that $q:L\to G$ is the kernel of $p:G\to A$. Then (a), (b)
follow from the well-know properties of normal subgroup (cf.
\cite[Chapter 13]{waterhause}).
 It remains to show (c).

Let ${\rm Ind} : {\rm Rep}(L) \to {\rm Rep}(G) $ be the induced
representation functor, it is the right adjoint functor to the
restriction functor ${\rm Res}:  {\rm Rep}(G)\to {\rm Rep}(L)$
that is, one has a functorial isomorphism 
\ga{a2}{\Hom_{G}(V,{\rm
Ind}( W))\xrightarrow{\cong} {\rm Hom}_L({\rm Res}( V),W) .} 
It is
easy to check 
\ga{a3}{\text{Ind}(W)\cong (k[G]\otimes_kV)^L} where
$L$ acts on $k[G]$ on the right. It is well-known that $k[G]$ is
faithfully flat over it subalgebra $k[A]$ (\cite{waterhause},
Chapter 13) and there is the following isomorphism
\ga{a4}{k[G]\otimes_{k[A]}k[G]\cong k[L]\otimes_kk[G]}
 which precisely means  that $G\to A$ is a principal bundle under $L$.
 Consequently
\ga{a5}{k[G]\otimes_{k[A]}\text{Ind}(W)\cong k[A]\otimes kV
 }
Thus the functor $\text{Ind}:\Rep(L)\to \Rep(G)$ is exact.

Setting $V={\rm Ind}( W)$ in \eqref{a2}, one obtains a canonical
map $u_W:{\rm Ind}( W)\to W$ in ${\rm Rep}(N)$ which gives back
the isomorphism in \eqref{a2} as follows: $\Hom_{G}(V,{\rm Ind}(
W))\ni h\mapsto u_W\circ h \in {\rm Hom}_N({\rm Res}( V),W)$. The
map $u_W$ is non-zero whenever $W$ is non-zero. Indeed, since
${\rm Ind}$ is faithfully exact, ${\rm Ind}(W)$ is non-zero
whenever $W$ is non-zero. Thus, if $u_W=0$, then  both sides of
\eqref{a2} are zero for any $V$. On the other hand,  for  $V={\rm
Ind}(W)$, the right hand side contains the identity map. This show
that $u_W$ can't vanish.

We want to show that $u_W$ is always surjective. Let $U={\rm
Im}(u_W)$ and $T=W/U$. We have the following diagram
\ga{a6}{\begin{CD}
0&@>>>& {\rm Ind}( U) &@>>> &{\rm Ind}( W) &@>>>&{\rm Ind}( T) &@>>>&0\\
&&&&@VVV&&@VVV&&@VVV\\
0&@>>>&U&@>>>&W&@>>>&T&@>>>&0
\end{CD}}
The composition ${\rm Ind}( W)\surj {\rm Ind}( T)\to T$ is 0,
therefore ${\rm Ind}( T)\to T$ is a zero map, implying $T=0$.

Since ${\rm Ind} (W)$ is a union of its finite dimensional
subrepresentations, we can therefore find a finite dimensional
$G$-subrepresentation $W_0(W)$ of ${\rm Ind} (W)$ which still maps
surjectively on $W$. In order to obtain the statement on the
embedding of $W$, we dualize $W_0(W^\vee)\surj W^\vee$.

Conversely, assume that (a), (b), (c) are satisfied. Then it
follows from (a) that for $U\in \Rep(A)$, $q^*p^*(U)\in \Rep(L)$
is trivial. Hence $pq:L\to A$ is the trivial homomorphism. Recall
that by assumption, $q$ is injective, $p$ is surjective. Let $\bar
q:\bar L\to G$ be the kernel of $p$. Then we have commutative
diagram 
\ga{a7}{\xymatrix{ L\ar[rr]^i \ar[rd]^q&&\bar
L\ar[dl]^{\bar q}\\ &G }\quad\Longrightarrow\quad \xymatrix{
\Rep(L)& &\Rep(\bar L)\ar[ll]^{i^*}\\ &\Rep(G)\ar[lu]^{q^*}\ar[ur]
^{\bar q^*} }} with injective homomorphisms. It remains to show
that $i$ is surjective, which amounts to saying that the category
$i^*\Rep(\bar L)$ in $\Rep(L)$ is full and closed under taking
subquotients.

We first show the fullness. Let $\bar W_0,\bar W_1$ be objects in
$\Rep(\bar L)$ and $\varphi:W_0:=i^*(\bar W_0)\to i^*(\bar
W_1)=:W_1$ be a morphism in $\Rep(L)$. Since $\Rep(\bar L)$ also
satisfies (c), there exists $V_0,\,V_1$ in $\Rep(G)$ with a
surjective morphism $\pi:\bar q^*(V_0)\to \bar W_0$, and an
injective morphism $\iota:\bar W_1\to \bar q^*(V_1)$. These yield
a morphism 
\ga{a8}{\psi:=i^*(\iota)\varphi i^*(\pi):q^*(V_0)\to
q^*(V_1)} The morphism $\psi$ induces and element in
$H^0(L,q^*(V_0^\vee\otimes V_1))$. Now, by (b) and by the fact
that $\Rep(\bar L)$
 also satisfies (b)
we conclude that $\psi=i^*(\bar \psi)$, for some $\bar\psi:\bar
q^*(V_0)\to \bar q^*(V_1)$. Since $\iota$ is injective and $\pi$
is surjective, we conclude that $\varphi=\bar\varphi$, for some
$\bar\varphi:\bar W_0\to \bar W_1$ in $\Rep(\bar L)$. Thus the
category $i^*\Rep(\bar L)$ is full in $\Rep(L)$.

On the other hand, for any $W\in \Rep(L)$, by (c) there exist
$V_0,V_1$ in $\Rep(G)$ and $\varphi:q^*(V_0)\to q^*(V_1)$ such
that $W\cong \im\varphi$. By the fullness of $i^*\Rep(\bar L)$ in
$\Rep(L)$, $\varphi=i^*\bar\varphi$, hence $W\cong
i^*(\im\bar\varphi)$. Thus we have proved that any object in
$\Rep(L)$ is isomorphic to the image under $i^*$ of an object in
$\Rep(\bar L)$. Together with the discussion above this implies
that $L\cong \bar L$.
\end{proof}
\end{appendix}

\bibliographystyle{plain}

\renewcommand\refname{References}

\end{document}